\documentclass[11pt]{article}
\usepackage[latin1]{inputenc}
\usepackage{amsmath,amsthm,amssymb}
\usepackage{amsfonts}
\usepackage{amsmath,amsthm,amssymb,amscd}
\usepackage{latexsym}
\usepackage{color}
\usepackage{graphicx}
\usepackage{mathrsfs}

\usepackage{color,enumitem,graphicx}
\usepackage[colorlinks=true,urlcolor=black,
citecolor=black,linkcolor=black,linktocpage,pdfpagelabels,
bookmarksnumbered,bookmarksopen]{hyperref}
\makeatletter \@addtoreset{equation}{section} \makeatother

\newtheorem{theorem}{Theorem}[section]

\newtheorem{proposition}{Proposition}[section]
\newtheorem{lemma}{Lemma}[section]
\newtheorem{remark}{Remark}[section]

\newtheorem{corollary}[theorem]{Corollary}

\allowdisplaybreaks

\begin{document}
\title{Normalized solutions for a Choquard equation with exponential growth in $\mathbb{R}^{2}$}

\author{Shengbing Deng\footnote{E-mail address:\, {\tt shbdeng@swu.edu.cn} (S. Deng), {\tt JwYumaths@163.com} (J. Yu)}\ \ \ and Junwei Yu\\
\footnotesize  School of Mathematics and Statistics, Southwest University,
Chongqing, 400715, P.R. China}

\date{ }
\maketitle

\begin{abstract}
{In this paper, we study the existence of normalized solutions to the following nonlinear Choquard equation with exponential growth
\begin{align*}
 \left\{
\begin{aligned}
&-\Delta u+\lambda u=(I_{\alpha}\ast F(u))f(u), \quad
\quad
\hbox{in }\mathbb{R}^{2},\\
&\int_{\mathbb{R}^{2}}|u|^{2}dx=a^{2},
\end{aligned}
\right.
\end{align*}
where $a>0$ is prescribed, $\lambda\in \mathbb{R}$, $\alpha\in(0,2)$, $I_{\alpha}$ denotes the Riesz potential, $\ast$ indicates the convolution operator, the function $f(t)$ has exponential growth in $\mathbb{R}^{2}$ and $F(t)=\int^{t}_{0}f(\tau)d\tau$. Using the Pohozaev manifold and variational methods, we establish the existence of normalized solutions to the above problem.}



\end{abstract}

\section{{\bfseries Introduction}}\label{introduction}

This paper concerns the existence of normalized solutions to the following nonlinear Choquard equation with exponential growth
\begin{equation}\label{aa}
 \left\{
\begin{aligned}
&-\Delta u+\lambda u=(I_{\alpha}\ast F(u))f(u), \quad
\quad
\hbox{in }\mathbb{R}^{2},\\
&\int_{\mathbb{R}^{2}}|u|^{2}dx=a^{2},
\end{aligned}
\right.
\end{equation}
where $a>0$ is prescribed, $\lambda\in \mathbb{R}$, $\alpha\in(0,2)$, $I_{\alpha}$ denotes the Riesz potential defined by
\begin{eqnarray*}
    \begin{aligned}\displaystyle
    I_{\alpha}(x)=\frac{\Gamma(\frac{2-\alpha}{2})}{\Gamma(\frac{\alpha}{2})2^{\alpha}\pi|x|^{2-\alpha}}:=\frac{A_{\alpha}}{|x|^{2-\alpha}}, \ \ x\in \mathbb{R}^{2}\backslash\{0\},
    \end{aligned}
\end{eqnarray*}
where $\Gamma$ represents the gamma function, $\ast$ indicates the convolution operator, $F(t)=\int^{t}_{0}f(\tau)d\tau$. The nonlinearity $f$ satisfies some suitable conditions that will be specified later.

Equation \eqref{aa} arises from seeking the standing wave solutions of prescribed mass to the following time-dependent nonlinear Choquard equation
\begin{equation}\label{abb}
\begin{aligned}
i \frac{\partial}{\partial t}\Phi=\Delta \Phi+(I_{\alpha}\ast F(\Phi))f(\Phi), \quad
\quad
\hbox{in }\mathbb{R}^{2}\times \mathbb{R},
\end{aligned}
\end{equation}
where $f(e^{i\lambda t}z)=e^{i\lambda t}f(z)$ and $F(e^{i\lambda t}z)=\int^{t}_{0}f(e^{i\lambda t}z)dz=e^{i\lambda t}F(z)$ for any $\lambda,t,z\in \mathbb{R} $. As we all know, an important feature of equation \eqref{abb} is conservation of mass with the $L^{2}$-norms $|\Phi(\cdot,t)|_{2}$ of solutions are independent of $t\in \mathbb{R}$. To find stationary states, one makes the ansatz $\Phi(x,t)=e^{i\lambda t}u(x)$, which still ensure conservation of mass with the $L^{2}$-norms.

We mention that Schrodinger equations with prescribed $L^{2}$-norms has received extensive attention in recent years. Jeanjean \cite{jeanjean1} considered the existence of solutions of prescribed mass to the following Schr\"{o}dinger equations
\begin{equation}\label{eqq1}
 \left\{
\begin{aligned}
&-\Delta u+\lambda u= f(u)  \quad \quad \hbox{in }\mathbb{R}^{N},\\
&\int_{\mathbb{R}^{N}}|u|^{2}dx=a^{2},
\end{aligned}
\right.
\end{equation}
where $a>0$ is prescribed, $N\geq1$ and $f$ satisfies some subcritical growth conditions. Using the minimax approach and a compactness argument, the author obtained the existence of normalized solutions for equations \eqref{eqq1}. In \cite{JeanjeanLus}, Jeanjean and Lu complemented and generalized the known results in \cite{jeanjean1}. In \cite{SOAVE1}, Soave studied the normalized solutions for equations \eqref{eqq1} with combined nonlinearities, where $f(t)=\mu|t|^{q-2}t+|t|^{p-2}t$ and $2<q\leq2+\frac{4}{N}\leq p <2^{*}$. For $f$ has a critical growth in the Sobolev sense, Soave \cite{SOAVE} obtained the existence of ground states and raised some questions with combined
nonlinearities. For more results in normalized solutions for the equations \eqref{eqq1}, the reader may refer to \cite{BartschSaovea,BartschSaoveb,JeanjeanLu,Lia,Wei,BartschDE,Ac} and references therein.

On the other hand, we mention that the study of our problem is based on some interesting results of the following Choquard equation
\begin{equation}\label{acc}
\begin{aligned}
&-\Delta u+V(x) u=(I_{\alpha}\ast F(u))f(u), \quad
\quad
\hbox{in }\mathbb{R}^{N}.
\end{aligned}
\end{equation}
This nonlocal equation plays an important role in quantum theory and description of finite range multi-body interaction. For $N=3,\alpha=2$ and $f(s)=s$, Fr$\ddot{o}$hlich \cite{FR} introduced this equation to study the modeling of quantum polaron. Moroz and Van Schaftingen \cite{MVA,MVa} obtained the existence of a ground state solution for equation \eqref{acc} and studied its some properties. For $N=2$ and $f$ has exponential critical growth, equation \eqref{acc} has been investigated by some authors; see e.g. \cite{ACTT,BATTA,QT} and the references therein. For more classical results regarding Choquard equation, we refer to \cite{MVb} for a good survey.

When the $L^{2}$-norms $|u|_{2}$ is prescribed, equation \eqref{acc} has important physical significance: in Bose-Einstein condensates and the nonlinear optics framework. In \cite{LiG}, Li and Ye firstly considered the following Choquard equation
\begin{equation}\label{c}
 \left\{
\begin{aligned}
&-\Delta u+\lambda u=(I_{\alpha}\ast F(u))f(u), \quad
\quad
\hbox{in }\mathbb{R}^{N},\\
&\int_{\mathbb{R}^{2}}|u|^{2}dx=a^{2},
\end{aligned}
\right.
\end{equation}
where $\lambda$ cannot be prescribed but appears as a Lagrange multiplier in a variational approach. The authors in \cite{LiG} used a minimax procedure and the concentration compactness to show equation \eqref{c} has at least a weak solution. For $N=3$, Yuan et al. \cite{Yuan} complemented and generalized the known results in \cite{LiG}. In \cite{BartschLiu}, Bartsch et al. proved the existence of a least energy solution of equation \eqref{c} in all dimensions $N\geq1$, which is simpler and more transparent than the one from \cite{LiG}. For critical case, Ye et al. \cite{Ye} obtained the existence of ground states for the critical Hartree equation with perturbation. For more results regarding equation \eqref{c}, the reader may refer to \cite{Lis,Yao} and references therein.

We know that classical Sobolev embedding that $H^{1}(\mathbb{R}^{N})$ is continuously embedded in $L^{q}(\mathbb{R}^{N})$ for all $q\in[2,2^{\ast}]$, where $2^{\ast}=2N/(N-2)$. Thus, we know that $2^{\ast}=\infty$, if $N=2$. In this case, Ji et al. \cite{AlvesJi} firstly studied equation \eqref{c} with exponential critical nonlinearities. Using the minimax approach, the authors firstly obtained results for normalized problem with two-dimensional exponential critical growth. In addition, when $N\geq3$, the authors in \cite{AlvesJi} complemented some recent results found in \cite{SOAVE}.

Now, let us introduce the precise assumptions under what our problem is studied. Assume $f$ satisfies the following conditions

($f_{1}$)$f(t):\mathbb{R}\rightarrow \mathbb{R}$ is continuous;

($f_{2}$) $f(t)=o(|t|^{\tau})$ as $|t|\rightarrow0$, for some $\tau>3$;

($f_{3}$) there exists a positive constant $\theta>2+\frac{\alpha}{2}$ such that $0<\theta F(t)\leq f(t)t$ for $t\neq0$;

($f_4$) there exist constants $\sigma>2+\frac{\alpha}{2}$ and $\mu>0$ such that
	$F(t)\geq \mu  \, |t|^{\sigma}$ for all $t >0$;

($f_5$) The function $\tilde{F}(t)=f(t)t-\frac{2+\alpha}{2}F(t)$ satisfies
	\begin{equation*}
		\frac{\tilde{F}(t)}{t^{2+\frac{\alpha}{2}}} \ \mbox{is nonincreasing in }(0,\infty) \ \mbox{and nondecreasing in }(0,\infty).
	\end{equation*}

In the subcritical case, we assume that $f$ satisfies the following condition

($f_6$) $f$ has exponential subcritical growth, i.e. for all $\gamma>0$ we have
\begin{eqnarray*}
    \begin{aligned}\displaystyle
    \lim\limits_{|t|\rightarrow +\infty}\frac{|f(t)}{e^{\gamma|t|^{2}}}=0.
    \end{aligned}
\end{eqnarray*}

The first result of this paper can be stated as follows:
\begin{theorem}\label{T1}
Assume that $f$ satisfies ($f_{1}$)-($f_{3}$) and ($f_{6}$). Then equation \eqref{aa} has a weak solution $(\lambda,u)$ with $\lambda>0$ and $u\in H^{1}(\mathbb{R}^{2})$. Moreover, if $(f_{5})$ is also assumed, then $u$ can be chosen as a nontrivial ground state solution of equation \eqref{aa}.
\end{theorem}

Motivated by the research made in the Choquard equations and \cite{AlvesJi}, in this paper we consider the exponential critical growth for Choquard equations in $\mathbb{R}^{2}$.  We recall that in $\mathbb{R}^{2}$, the natural growth restriction on the function $f$ is given by the Trudinger-Moser inequality.

In the critical case, we assume that $f$ satisfies the following condition

($f_7$) $f$ has $\gamma_{0}$-exponential critical growth, i.e. there exists $\gamma_{0}>0$ such that
\begin{align*}
    \lim\limits_{|t|\rightarrow +\infty}\frac{|f(t)|}{e^{\gamma|t|^{2}}}=
 \left\{
    \begin{aligned}
    0, \ \ \ \forall \gamma> \gamma_{0},\\
    +\infty, \ \ \ \forall \gamma< \gamma_{0}.
    \end{aligned}
\right.
\end{align*}

Our second result is as follows:

\begin{theorem}\label{T2}
Assume that $f$ satisfies ($f_{1}$)-($f_{4}$) and ($f_7$). If $a^{2} <\frac{(2+\alpha)\pi}{\gamma_{0}}$, then there exists $\mu^*=\mu^*(a)>0$ such that equation \eqref{aa} has a weak solution $(\lambda,u)$ with $\lambda>0$ and $u\in H^{1}(\mathbb{R}^{2})$ for all $\mu \geq \mu^*$. Moreover, if $(f_{5})$ is also assumed, then $u$ can be chosen as a nontrivial ground state solution of equation \eqref{aa}.
\end{theorem}

\begin{remark}
There are many functions satisfying $(f_{1})$-$(f_{5})$ and $(f_{7})$. See, for example
\begin{eqnarray*}
    \begin{aligned}\displaystyle
    f(t)=|t|^{\sigma-2}t+|t|^{\tau}(e^{t^{2}}-1) \ \ \mbox{for} \ \ t\in \mathbb{R},
    \end{aligned}
\end{eqnarray*}
where $\sigma>2+\frac{\alpha}{2}$ and $\tau>3$.
\end{remark}

\begin{remark}
In this article, we suppose that $f$ satisfies $(f_{6})$ or $(f_{7})$, which means $f$ has subcritical exponential growth or critical exponential growth. This notion of criticality is motivated by the inequality of Trudinger and Moser; see \cite{Moser,Trudinger}.
\end{remark}

The paper is organized as follows: In Section $\ref{preliminaries}$, the variational setting and some preliminary results are presented. In Section $\ref{Sgeo}$, we introduce the geometry structure related to equation \eqref{aa}. In Section $\ref{PSa}$, we give some properties about the $(PS)$ sequence. In Section $\ref{PFTa}$, we complete the proof of subcritical case. Section $\ref{minimax}$ is devoted to give an estimate for the minimax level of the critical case. Finally in Section $\ref{PFTb}$, we complete the proof of critical case.

\vspace{0.5 cm}

\noindent \textbf{Notation:} From now on in this paper, otherwise mentioned, we use the following notations:
\begin{itemize}
    \item$C,C_1,C_2,...$ denote some positive constants.
	\item  $|\,\,\,|_p$ denotes the usual norm of the Lebesgue space $L^{p}(\mathbb{R}^2)$, for $p \in [1,+\infty]$.
	\item $\Vert\,\,\,\Vert$ denotes the usual norm of the Sobolev space $H^{1}(\mathbb{R}^2)$.
    \item $o_{n}(1)$ denotes a real sequence with $o_{n}(1)\to 0$ as $n \rightarrow+\infty$.
\end{itemize}

\section{{\bfseries Preliminaries and functional setting}}\label{preliminaries}
In this section, we give some preliminary results and outline the variational framework for \eqref{aa}.
\begin{proposition}\label{PRa}
(\cite{Cao}).$i)$ If $\gamma>0$ and $u\in H^{1}(\mathbb{R}^{2})$, then
\begin{eqnarray*}
    \begin{aligned}\displaystyle
    \int_{\mathbb{R}^{2}}\Big(e^{\gamma |u|^{2}}-1\Big)dx<\infty;
    \end{aligned}
\end{eqnarray*}
$ii)$ if $u\in H^{1}(\mathbb{R}^{2})$, $|\nabla u|^{2}_{2}\leq1$, $|u|_{2}\leq M<\infty$, and $\gamma<4\pi$, then there exists a constant $\mathcal{C}(M,\gamma)$, which depends only on $M$ and $\alpha$, such that
\begin{eqnarray*}
    \begin{aligned}\displaystyle
    \int_{\mathbb{R}^{2}}\Big(e^{\gamma |u|^{2}}-1\Big)dx<\mathcal{C}(M,\gamma).
    \end{aligned}
\end{eqnarray*}
\end{proposition}

Now we recall the Hardy-Littlewood-Sobolev inequality, see \cite{Lieb}.

\begin{proposition}\label{PRc}(Hardy-Littlewood-Sobolev inequality). Let $t,r>1$, $0<\alpha<2$, with $\frac{1}{t}+\frac{2-\alpha}{2}+\frac{1}{r}=2$, $f\in L^{t}(\mathbb{R}^{2})$ and $g\in L^{r}(\mathbb{R}^{2})$. Then there exists a sharp constant $C(t,\alpha,r)$, independent of $f,g$, such that
\begin{eqnarray*}
    \begin{aligned}\displaystyle
    \int_{\mathbb{R}^{2}} (I_{\alpha}\ast f)g dx \leq C(t,\alpha,r)|f|_{t}|g|_{r}.
    \end{aligned}
\end{eqnarray*}
\end{proposition}

Applying Proposition \ref{PRc}, we know
\begin{eqnarray*}
    \begin{aligned}\displaystyle
    \int_{\mathbb{R}^{2}} (I_{\alpha}\ast F(u))F(u) dx
    \end{aligned}
\end{eqnarray*}
is well defined if $F(u)\in L^{t}(\mathbb{R}^{2})$ for $t>1$ given by
\begin{eqnarray*}
    \begin{aligned}\displaystyle
    \frac{2}{t}+\frac{2-\alpha}{2}=2.
    \end{aligned}
\end{eqnarray*}
This implies that we must require
\begin{eqnarray*}
    \begin{aligned}\displaystyle
    F(u)\in L^{\frac{4}{2+\alpha}}(\mathbb{R}^{2}).
    \end{aligned}
\end{eqnarray*}

In order to apply variational methods, we recall that $H^{1}(\mathbb{R}^{2})$ denotes the usual Sobolev space with the inner product and norm
\begin{eqnarray*}
    \begin{aligned}\displaystyle
    (u,v)=\int_{\mathbb{R}^{2}}(\nabla u\nabla v +uv)dx, \ \ \|u\|=( u,u)^{1/2}, \ \ \forall \ u,v \in H^{1}(\mathbb{R}^{2}).
    \end{aligned}
\end{eqnarray*}

Solutions of equation \eqref{aa} correspond to critical points of the energy functional $J:H^{1}(\mathbb{R}^{2})\rightarrow \mathbb{R}$ defined by
\begin{equation}\label{func}
    \begin{aligned}\displaystyle
J(u)=\frac{1}{2}\int_{\mathbb{R}^{2}} |\nabla u|^2 dx-\frac{1}{2}\int_{\mathbb{R}^{2}} (I_{\alpha}\ast F(u))F(u)dx
\end{aligned}
\end{equation}
and constrained to the $L^{2}$-torus
\begin{eqnarray*}
    \begin{aligned}\displaystyle
    \mathcal{T}(a):=\{u\in H^{1}(\mathbb{R}^{2}): |u|_{2}=a\}.
    \end{aligned}
\end{eqnarray*}
The parameters $\lambda\in \mathbb{R}$ will appear as Lagrange multiplier. By Proposition \ref{PRa} and \ref{PRc}, we know that $(I_{\alpha}\ast F(u))F(u)\in L^{1}(\mathbb{R}^{2})$, which implies that $J$ is well defined. It is easy to see that $J$ is of class $C^{1}$, and that it is unbounded from below on $\mathcal{T}(a)$. It is well known that critical points of $J$ will not satisfy the Palais-Smale condition, as a consequence we recall that solutions of \eqref{aa} satisfy the Pohozaev identity
\begin{equation}\label{Pohozaev}
    \begin{aligned}\displaystyle
    P(u):=\int_{\mathbb{R}^{2}} |\nabla u|^2 dx+\frac{2+\alpha}{2}\int_{\mathbb{R}^{2}} (I_{\alpha}\ast F(u))F(u) dx -\int_{\mathbb{R}^{2}} (I_{\alpha}\ast F(u))f(u)u dx.
    \end{aligned}
\end{equation}

Now, we introduce the $L^{2}$-invariant scaling $s\star u(x):=e^{s}u(e^{s}x)$. To be more precise, for $u\in \mathcal{T} (a)$ and $s\in \mathbb{R}$, let $\mathcal{H}: H^{1}(\mathbb{R}^{2})\times \mathbb{R}\rightarrow H^{1}(\mathbb{R}^{2})$ defined by
\begin{eqnarray*}
    \begin{aligned}\displaystyle
    \mathcal{H}(u, s)(x)=e^{s}u(e^{s}x).
    \end{aligned}
\end{eqnarray*}
Then $\tilde{J}: \mathbb{R} \rightarrow \mathbb{R}$ with
\begin{eqnarray*}
    \begin{aligned}\displaystyle
	\tilde{J}_{u}(s):=J(\mathcal{H}(u, s))=\frac{e^{2s}}{2}\int_{\mathbb{R}^{2}} |\nabla u|^2 dx-\frac{1}{2e^{(2+\alpha)s}}\int_{\mathbb{R}^{2}} (I_{\alpha}\ast F(e^{s}u(x)))F(e^{s}u(x))\,dx.
    \end{aligned}
\end{eqnarray*}
A straightforward calculation shows that $\tilde{J}^{\prime}_{u}(s)=P(\mathcal{H}(u, s))$.

In order to overcome the loss of compactness of the Sobolev embedding in whole $\mathbb{R}^{2}$, in our opinion, we work on the space $H^{1}_{rad}(\mathbb{R}^{2})$ to get some compactness results. Thus, we define
\begin{eqnarray*}
    \begin{aligned}\displaystyle
     \mathcal{T}_{r}(a):=\mathcal{T}(a)\cap H^{1}_{rad}(\mathbb{R}^{2}).
    \end{aligned}
\end{eqnarray*}

From $(f_{1})$-$(f_{3})$ and if $f(t)$ has subcritical exponential growth at $+\infty$, we have the following immediate result: fix $q>2$ and $\tau>3$, for any $\varepsilon>0$ and $\gamma>0$, there exists a constant $\kappa_{\varepsilon}>0$, which depends on $q,\gamma,\varepsilon$, such that
\begin{equation}\label{snonb}
    \begin{aligned}\displaystyle
| f(t)|\leq\varepsilon|t|^{\tau}+\kappa_{\varepsilon}|t|^{q-1}(e^{\gamma |t|^{2}}-1) \ \mbox{for all} \ t\in \mathbb{R},
    \end{aligned}
\end{equation}
and
\begin{equation}\label{snona}
    \begin{aligned}\displaystyle
|F(t)|\leq\varepsilon|t|^{\tau+1}+\kappa_{\varepsilon}|t|^{q}(e^{\gamma |t|^{2}}-1) \ \mbox{for all} \ t\in \mathbb{R}.
    \end{aligned}
\end{equation}
Similarly, if $f(t)$ has critical exponential growth at $+\infty$ with critical exponent $\gamma_{0}$, then fix $q>2$ and $\tau>3$, for any $\varepsilon>0$ and $\gamma>\gamma_{0}$, there exists a constant $\kappa_{\varepsilon}>0$, which depends on $q,\gamma,\varepsilon$ and $\mu$, such that
\begin{equation}\label{nonb}
    \begin{aligned}\displaystyle
| f(t)|\leq\varepsilon|t|^{\tau}+\kappa_{\varepsilon}|t|^{q-1}(e^{\gamma |t|^{2}}-1) \ \mbox{for all} \ t\in \mathbb{R},
    \end{aligned}
\end{equation}
and
\begin{equation}\label{nona}
    \begin{aligned}\displaystyle
|F(t)|\leq\varepsilon|t|^{\tau+1}+\kappa_{\varepsilon}|t|^{q}(e^{\gamma |t|^{2}}-1) \ \mbox{for all} \ t\in \mathbb{R}.
    \end{aligned}
\end{equation}

In Section \ref{Sgeo} and \ref{PSa}, we consider equation \eqref{aa} in subcritical case or critical case. By $(f_{6})$, we know that $\gamma \| u\|^{2}  < (2+\alpha)\pi$ for $\gamma>0$ close to $0$ and $u$ is bounded in $H^{1}(\mathbb{R}^{2})$. Similarly, by $(f_{7})$, we know that $\gamma \| u\|^{2}  < (2+\alpha)\pi$ for $\gamma>\gamma_{0}$ close to $\gamma_{0}$ and $\|u\|^{2}<\frac{(2+\alpha)\pi}{\gamma_{0}}$. Hence, in what follows, we consider the problem for a unified condition $\gamma \| u\|^{2}  < (2+\alpha)\pi$.

\section{{\bfseries The minimax approach}}\label{Sgeo}
To find a solution of \eqref{aa}, in this section, we show that $\tilde{J}$ on $\mathcal{T}_{r}(a)\times \mathbb{R}$ possesses a kind of mountain-pass geometrical structure.
\begin{lemma}\label{geo} Assume that $(f_{1})$-$(f_{3})$ hold. Let $u\in \mathcal{T}_{r}(a)$. Then we have:

(i)$|\nabla\mathcal{H}(u, s)|_{2}\rightarrow 0$ and $J(\mathcal{H}(u, s))\rightarrow 0$ as $s\rightarrow -\infty$;

(ii)$|\nabla\mathcal{H}(u, s)|_{2}\rightarrow +\infty$ and $J(\mathcal{H}(u, s))\rightarrow -\infty$ as $s\rightarrow +\infty$.
\end{lemma}

\begin{proof} \mbox{} By a straightforward calculation, it follows that
\begin{equation}\label{cria}
    \int_{\mathbb{R}^2}|e^{s}u(e^{s}x)|^{2}\,dx=a^{2}, \quad
    \int_{\mathbb{R}^2} |\nabla\mathcal{H}(u,s)(x)|^{2}\,dx=e^{2s}\int_{\mathbb{R}^2}|\nabla u|^2dx,
\end{equation}
and
\begin{equation}\label{3.2}
\int_{\mathbb{R}^2}|\mathcal{H}(u, s)|^{\xi}\,dx= e^{(\xi-2)s}\int_{\mathbb{R}^2}|u(x)|^{\xi}\,dx, \quad \forall \xi \geq 2.
\end{equation}

From the above equalities, fixing $\xi>2$, we have
\begin{equation}\label{3.4}
    |\nabla\mathcal{H}(w,s)|_{2}\rightarrow 0 \quad \mbox{and} \quad |\mathcal{H}(w,s)|_{\xi} \rightarrow 0 \quad  \mbox{as} \quad s \to -\infty.
	\end{equation}
By \eqref{snona} and \eqref{nona}, we obtain
\begin{eqnarray*}
    \begin{aligned}\displaystyle
|F(u)|\leq\varepsilon |u|^{\tau+1}+\kappa_{\varepsilon}|u|^{q}(e^{\gamma |u|^{2}}-1)\,\, \text{ for all }\, u \in H^{1}(\mathbb{R}^{2}),
    \end{aligned}
\end{eqnarray*}
where $q>2$. Hence,
\begin{eqnarray*}
    \begin{aligned}\displaystyle
|F(\mathcal{H}(u,s))| \leq \varepsilon| \mathcal{H}(u, s)|^{\tau+1}+\kappa_{\varepsilon}| \mathcal{H}(u, s)|^{q}(e^{\gamma | \mathcal{H}(u, s)|^{2}}-1).
    \end{aligned}
\end{eqnarray*}
By Proposition \ref{PRa}, one has
\begin{eqnarray*}
    \begin{aligned}\displaystyle
    \int\limits_{\mathbb{R}^{2}}(e^{\gamma|\mathcal{H}(u,s)|^{2}}-1)dx&=\int\limits_{\mathbb{R}^{2}}(e^{\gamma\|\mathcal{H}(u, s)\|^{2}(\mathcal{H}(u, s)/\|\mathcal{H}(u, s)\|)^{2}}-1)dx\leq C.
    \end{aligned}
\end{eqnarray*}
for all $\gamma\|\mathcal{H}(u, s)\|^{2}< (2+\alpha)\pi$. Hence, using Proposition \ref{PRc} and H\"{o}lder equality, we deduce that
\begin{eqnarray*}
    \begin{aligned}\displaystyle
\Big(\int_{\mathbb{R}^{2}}|F(\mathcal{H}(u,s))|^{\frac{4}{2+\alpha}}dx\Big)^{\frac{2+\alpha}{4}}
&\leq \Big(\int_{\mathbb{R}^{2}}|\varepsilon|\mathcal{H}(u,s)|^{\tau+1}+C|\mathcal{H}(u,s)|^{q}(e^{\gamma |\mathcal{H}(u,s)|^{2}}-1)|^{\frac{4}{2+\alpha}}dx\Big)^{\frac{2+\alpha}{4}}\\
&\leq  |\varepsilon|\mathcal{H}(u,s)|^{\tau+1}|_{\frac{4}{2+\alpha}}+|C|\mathcal{H}(u,s)|^{q}(e^{\gamma |\mathcal{H}(u,s)|^{2}}-1)|_{\frac{4}{2+\alpha}}   \\
&\leq  \varepsilon|\mathcal{H}(u,s)|^{\tau+1}_{\frac{4(\tau+1)}{2+\alpha}}+C|\mathcal{H}(u,s)|^{q}_{\frac{4qt^{\prime}}{2+\alpha}}(\int_{\mathbb{R}^{2}}(e^{\frac{4}{2+\alpha}\gamma t|\mathcal{H}(u,s)|^{2}}-1)dx)^{\frac{2+\alpha}{4t}}.
    \end{aligned}
\end{eqnarray*}
Then, there exists $t>1$ close to $1$ such that
\begin{eqnarray*}
    \begin{aligned}\displaystyle
	\gamma t\|\mathcal{H}(u, s)\|^{2}\leq (2+\alpha)\pi,
    \end{aligned}
\end{eqnarray*}
which implies that
\begin{equation} \label{Domina1}
    \begin{aligned}\displaystyle
\Big(\int_{\mathbb{R}^{2}}|F(\mathcal{H}(u,s))|^{\frac{4}{2+\alpha}}dx\Big)^{\frac{2+\alpha}{4}}\leq C.
    \end{aligned}
\end{equation}
Then, by $\eqref{Domina1}$, we have
\begin{eqnarray*}
    \begin{aligned}\displaystyle
\int_{\mathbb{R}^2}(I_{\alpha}\ast F(\mathcal{H}(u, s)))F(\mathcal{H}(u, s))dx &\leq |F(\mathcal{H}(u, s))|_{\frac{4}{2+\alpha}}|F(\mathcal{H}(u, s))|_{\frac{4}{2+\alpha}}\\
&\leq \Big( \varepsilon\int_{\mathbb{R}^2}| \mathcal{H}(u, s)|^{\tau+1}dx+C(\int_{\mathbb{R}^2}| \mathcal{H}(u,s)|^{qt'}dx)^{1/t'} \Big)^{2},
    \end{aligned}
\end{eqnarray*}
where $t'=\frac{t}{t-1}$, and $t>1$ close to $1$.  Now, by using (\ref{3.2}), we obtain
\begin{eqnarray*}
    \begin{aligned}\displaystyle
\int_{\mathbb{R}^2}(I_{\alpha}\ast F(\mathcal{H}(u, s)))F(\mathcal{H}(u, s))dx \leq \Big( \varepsilon e^{(\tau-1)s}\int_{\mathbb{R}^2} |u(x)|^{\tau+1}dx+C(e^{(qt^{\prime}-2)s}\int_{\mathbb{R}^2}|u(x)|^{qt'}dx)^{1/t'}\Big)^{2}.
    \end{aligned}
\end{eqnarray*}
Thus, by $\tau-1>0$ and $qt^{\prime}-2>0$, we know that
\begin{eqnarray*}
    \begin{aligned}\displaystyle
	\int_{\mathbb{R}^2}(I_{\alpha}\ast F(\mathcal{H}(u, s)))F(\mathcal{H}(u, s))dx \rightarrow 0 \quad \mbox{as} \quad s \to -\infty,
    \end{aligned}
\end{eqnarray*}
which implies that $J(\mathcal{H}(u, s))\rightarrow 0$ as $s\rightarrow -\infty$, showing $(i)$.

In order to show $(ii)$, we define
\begin{eqnarray*}
    \begin{aligned}\displaystyle
 g(z)= \int_{\mathbb{R}^{2}}(I_{\alpha}\ast F(z))F(z)dx.
    \end{aligned}
\end{eqnarray*}
Then, note that by (\ref{cria}),
\begin{eqnarray*}
    \begin{aligned}\displaystyle
	J(\mathcal{H}(u, s))=\frac{e^{2s}}{2}\int_{\mathbb{R}^{2}} |\nabla u|^2 dx-\frac{1}{e^{(2+\alpha)s}}\int_{\mathbb{R}^{2}} (I_{\alpha}\ast F(e^{s}u(x)))F(e^{s}u(x))\,dx.
    \end{aligned}
\end{eqnarray*}
Set
\begin{eqnarray*}
    \begin{aligned}\displaystyle
w(t)=g(\frac{tu}{\|u\|}),
    \end{aligned}
\end{eqnarray*}
where $t=e^{s}$. By ($f_{3}$) we know
\begin{eqnarray*}
    \begin{aligned}\displaystyle
\frac{w^{\prime}(t)}{w(t)}\geq\frac{2\theta}{t} \ \mbox{for} \ t>0,
    \end{aligned}
\end{eqnarray*}
which implies that
\begin{eqnarray*}
    \begin{aligned}\displaystyle
g(tu)\geq g(\frac{u}{\|u\|})\|u\|^{2\theta}t^{2\theta}.
    \end{aligned}
\end{eqnarray*}
Therefore, we obtain
\begin{eqnarray*}
    \begin{aligned}\displaystyle
	J(\mathcal{H}(u, s)) \leq C_{1}e^{2s}-C_{2}e^{(2\theta-(2+\alpha))s}.
    \end{aligned}
\end{eqnarray*}
Since $2\theta-(2+\alpha)>2$, the last inequality yields $J(\mathcal{H}(w, s))\rightarrow -\infty$ as $s\rightarrow +\infty$.
\end{proof}

\begin{lemma}  \label{P1} Assume that $(f_{1})$-$(f_{3})$ hold. Let $u\in \mathcal{T}_{r}(a)$. Then there exists $K(a)>0$ small enough such that
\begin{eqnarray*}
    \begin{aligned}\displaystyle
	0<\sup_{u\in A} J(u)<\inf_{u\in B} J(u)
    \end{aligned}
\end{eqnarray*}
	with
\begin{eqnarray*}
    \begin{aligned}\displaystyle
	A=\left\{u\in \mathcal{T}_{r}(a), \int_{\mathbb{R}^2} |\nabla u|^2 dx\leq K(a) \right\}\quad \mbox{and} \quad B=\left\{u\in \mathcal{T}_{r}(a), \int_{\mathbb{R}^2} |\nabla u|^2 dx=2K(a) \right\}.
    \end{aligned}
\end{eqnarray*}
\end{lemma}

\begin{proof}
We will need the following Gagliardo-Sobolev inequality: for any $p> 2$,
\begin{eqnarray*}
    \begin{aligned}\displaystyle
	|u|_{p}\leq C(p, 2)|\nabla u|_2^{d_{p}}|u|_2^{1-d_{p}},
    \end{aligned}
\end{eqnarray*}
where $d_{p}=2(\frac{1}{2}-\frac{1}{p})$ and $u\in H^{1}(\mathbb{R}^{2})$. Now, let $K(a)<\frac{(2+\alpha)\pi}{2\gamma}-\frac{a^{2}}{2}$. Assume that $u_{1},u_{2}\in \mathcal{T}_{r}(a)$ such that $|\nabla u_{1}|^{2}_{2}\leq K(a)$ and $|\nabla u_{2}|^{2}_2=2K(a)$, Thus, we have
\begin{eqnarray*}
    \begin{aligned}\displaystyle
\|u_{2}\|=|\nabla u_{2}|^{2}_{2}+| u_{2}|^{2}_{2}<\frac{(2+\alpha)\pi}{\gamma}.
    \end{aligned}
\end{eqnarray*}
Then, similar as Lemma \ref{geo}, we obtain
\begin{eqnarray*}
    \begin{aligned}\displaystyle
	\int_{\mathbb{R}^2}(I_{\alpha}\ast F(u))F(u)\,dx \leq \varepsilon|u|^{\tau+1}_{\tau+1}+C_2|u|^{q}_{qt'},
    \end{aligned}
\end{eqnarray*}
where $\tau>3$, $q>2$, $t'=\frac{t}{t-1}>1$. By the Gagliardo-Sobolev inequality, we have
\begin{eqnarray*}
    \begin{aligned}\displaystyle
	\int_{\mathbb{R}^2}(I_{\alpha}\ast F(u_{2}))F(u_{2})\,dx \leq C_{1}K(a)^{\frac{\tau-1}{2}}+C_2K(a)^{(\frac{q}{2}-\frac{1}{t'})}.
    \end{aligned}
\end{eqnarray*}
From $(f_{3})$, we have $(I_{\alpha}\ast F(u_{1}))F(u_{1})> 0$ for any $u_{1}\in H^{1}(\mathbb{R}^{2})$. Then,
\begin{eqnarray*}
    \begin{aligned}\displaystyle
J(u_{2})-J(u_{1})&\geq\frac{1}{2}\int_{\mathbb{R}^2}|\nabla u_{2}|^{2}\,dx-\frac{1}{2}\int_{\mathbb{R}^2}|\nabla u_{1}|^{2}\,dx-\int_{\mathbb{R}^2}(I_{\alpha}\ast F(u_{2}))F(u_{2})\,dx\\
&\geq K(a)-\frac{1}{2}K(a) -C_{3}K(a)^{\frac{\tau-1}{2}}-C_{4} K(a)^{(\frac{q}{2}-\frac{1}{t'})}.
    \end{aligned}
\end{eqnarray*}
Hence,
\begin{eqnarray*}
    \begin{aligned}\displaystyle
	J(u_{2})-J(u_{1}) \geq \frac{1}{2}K(a)-C_{3}K(a)^{\frac{\tau-1}{2}}-C_{4} K(a)^{(\frac{q}{2}-\frac{1}{t'})}.
    \end{aligned}
\end{eqnarray*}
Since $\tau>3$ and $t'>0$ with $\frac{q}{2}-\frac{1}{t'}>1$, choosing $K(a)$ small enough if necessary, it follows that
\begin{eqnarray*}
    \begin{aligned}\displaystyle
	\frac{1}{2}K(a)-C_3K(a)^{\frac{\tau-1}{2}}-C_4 K(a)^{(\frac{q}{2}-\frac{1}{t'})}>0,
    \end{aligned}
\end{eqnarray*}
	showing the desired result.
\end{proof}

As a byproduct of the last lemma is the following corollary.
\begin{corollary} \label{newcor}
Assume that $(f_{1})$-$(f_{3})$ hold. Let $u\in \mathcal{T}_{r}(a)$. Then, if $|\nabla u|^{2}_{2}\leq K(a)$, there exists $K(a)>0$ small enough such that $J(u) >0$. Moreover,
\begin{eqnarray*}
    \begin{aligned}\displaystyle
	J_{\ast}=\inf\Big\{J(u): u\in \mathcal{T}_{r}(a) \ \ \mbox{and} \ \ |\nabla u|^{2}_{2}= \frac{K(a)}{2}\Big\}>0.
    \end{aligned}
\end{eqnarray*}
\end{corollary}

\begin{proof} Arguing as in the last lemma, we have for all $w\in \mathcal{T}_{r}(a)$
\begin{eqnarray*}
    \begin{aligned}\displaystyle
	J(u) \geq \frac{1}{2}|\nabla u|_{2}^{2}-C_{1}|\nabla u|_2^{\tau-1}-C_{2}|\nabla u|^{(q-\frac{2}{t^{\prime}})}_{2}>0,
    \end{aligned}
\end{eqnarray*}
	for $K(a)>0$ small enough.	
\end{proof}

In what follows, we fix $u_{0} \in \mathcal{T}_{r}(a)$ and apply Lemma \ref{geo}, \ref{P1} and Corollary  \ref{newcor}  to get two numbers $s_{1}<0$ and $s_{2}>0$, of such way that the functions $u_{1}=\mathcal{H}(u_{0}, s_{1})$ and $u_{2}=\mathcal{H}(u_{0}, s_{2})$ satisfy
\begin{eqnarray*}
    \begin{aligned}\displaystyle
|\nabla u_{1}|^2_2<\frac{K(a)}{2},\,\,  |\nabla u_{2}|_{2}^{2}>2K(a),\,\, J(u_{1})>0\,\,\mbox{and} \,\, J(u_{2})<0.
    \end{aligned}
\end{eqnarray*}

Now, following the ideas from Jeanjean \cite{jeanjean1}, we fix the following mountain pass level given by
\begin{eqnarray*}
    \begin{aligned}\displaystyle
m(a)=\inf_{h \in \Gamma}\max_{t \in [0,1]}J(h(t)),
    \end{aligned}
\end{eqnarray*}
where
\begin{eqnarray*}
    \begin{aligned}\displaystyle
\Gamma=\left\{h \in C([0,1], \ \mathcal{T}_{r}(a))\,:\,h(0)=u_{1} \,\, \mbox{and} \,\, h(1)=u_{2} \right\}.
    \end{aligned}
\end{eqnarray*}
By Corollary \ref{newcor}, we have
\begin{eqnarray*}
    \begin{aligned}\displaystyle
\max_{t \in [0,1]}J(h(t))\geq J_{\ast}>0.
    \end{aligned}
\end{eqnarray*}
Then, we obtain that $m(a)\geq J_{\ast}>0$.

\section{{\bfseries Palais-Smale sequence}}\label{PSa}
In this section, we take $\{u_{n}\}\subset \mathcal{T}_{r}(a)$ demotes the $(PS)$ sequence associated with the level $m(a)$ for $J$. Using $\hat{u}_{n}=\mathcal{H}(u_{n},s_{n})$, we know that $\{\hat{u}_{n}\}$ is a $(PS)$ sequence associated with the level $m(a)$ for $\tilde{J}$. Thus, we have
\begin{equation} \label{mu2}
	J(u_{n}) \to m(a), \,\, \mbox{as} \,\, n \to +\infty,
\end{equation}
\begin{equation} \label{EQ2}
	-\Delta u_{n}+\lambda_{n}u_{n}=(I_{\alpha}\ast F(u_{n}))f(u_{n})\ + o_n(1),
\end{equation}
for some sequence $\{\lambda_{n}\} \subset \mathbb{R}$, and
\begin{equation} \label{PEQ2}
	P(u_{n})=\int_{\mathbb{R}^2}|\nabla u_{n}|^{2} dx +\frac{2+\alpha}{2}\int_{\mathbb{R}^{2}} (I_{\alpha}\ast F(u_{n}))F(u_{n}) dx- \int_{\mathbb{R}^2} (I_{\alpha}\ast F(u_{n}))f(u_{n})u_{n},
\end{equation}
where $P(u_{n})\to 0$ as $n \to +\infty$. Then, setting
\begin{eqnarray*}
    \begin{aligned}\displaystyle
    \mathcal{P}(a):=\{u\in \mathcal{T}_{r}(a):P(u)=0\},
    \end{aligned}
\end{eqnarray*}
we are interested in the problem whether
\begin{eqnarray*}
    \begin{aligned}\displaystyle
    c(a):=\inf\limits_{u\in \mathcal{P}(a)}J(u)
    \end{aligned}
\end{eqnarray*}
is achieved. Then, in order to prove it, we give some properties of the $(PS)$ sequence.

\begin{lemma}\label{bounded}Assume that ($f_{1}$)-($f_{3}$) hold. Let $u\in \mathcal{T}_{r}(a)$. Then the $(PS)$ sequence $\{u_{n}\}$ of $J$ is bounded in $H^{1}(\mathbb{R}^{2})$.
\end{lemma}

\begin{proof}
Combining \eqref{mu2} and \eqref{PEQ2}, we obtain
\begin{eqnarray*}
    \begin{aligned}\displaystyle
(2+\alpha)m(a)+o_{n}(1)= (2+\alpha)J(u_{n})+P(u_{n})=(2+\frac{\alpha}{2})|\nabla u_{n}|_{2}^{2}-\int_{\mathbb{R}^2} (I_{\alpha}\ast F(u_{n}))f(u_{n})u_{n} dx.
    \end{aligned}
\end{eqnarray*}
By ($f_{3}$) we have
\begin{eqnarray*}
    \begin{aligned}\displaystyle
    \int_{\mathbb{R}^2} (I_{\alpha}\ast F(u_{n}))f(u_{n})u_{n} dx\geq \theta \int_{\mathbb{R}^{2}} (I_{\alpha}\ast F(u_{n}))F(u_{n}) dx.
    \end{aligned}
\end{eqnarray*}
Thus we obtain
\begin{equation}\label{bound1}
    \begin{aligned}\displaystyle
(2+\frac{\alpha}{2})|\nabla u_{n}|_{2}^{2}-\theta \int_{\mathbb{R}^{2}} (I_{\alpha}\ast F(u_{n}))F(u_{n}) dx \geq (2+\alpha)m(a).
    \end{aligned}
\end{equation}
On the other hand, using again \eqref{mu2}, it follows that
\begin{equation}\label{bound2}
    \begin{aligned}\displaystyle
    |\nabla u_{n}|_{2}^{2}= 2m(a)+\int_{\mathbb{R}^{2}} (I_{\alpha}\ast F(u_{n}))F(u_{n}) dx.
    \end{aligned}
\end{equation}
Taking \eqref{bound1} in \eqref{bound2}, we deduce that
\begin{eqnarray*}
    \begin{aligned}\displaystyle
    (2+\frac{\alpha}{2}-\theta)\int_{\mathbb{R}^{2}} (I_{\alpha}\ast F(u_{n}))F(u_{n}) dx\geq -2m(a).
    \end{aligned}
\end{eqnarray*}
Now, owing to $\theta>2+\frac{\alpha}{2}$, we get $\int_{\mathbb{R}^{2}} (I_{\alpha}\ast F(u_{n}))F(u_{n}) dx$ is bounded. Then, by \eqref{bound2}, $|\nabla u_{n}|_{2}^{2}$ is bounded. From the definition of norm in $H^{1}(\mathbb{R}^{2})$, we have
\begin{eqnarray*}
    \begin{aligned}\displaystyle
    \|u_{n}\|^{2}=|\nabla u_{n}|_{2}^{2}+|u_{n}|_{2}^{2}=|\nabla u_{n}|_{2}^{2}+a^{2}.
    \end{aligned}
\end{eqnarray*}
Then, the $(PS)$ sequence $\{u_{n}\}$ is bounded in $H^{1}(\mathbb{R}^{2})$.
\end{proof}

By Proposition \ref{PRa}, similar to \cite[Lemma 4.1]{AlvesJi}, we have the following Lemma.
\begin{lemma}\label{TMineq} Let $\{u_{n}\}$ be a sequence in $H^{1}(\mathbb{R}^{2})$ with $u_{n}\in \mathcal{T}_{r}(a)$ and
\begin{eqnarray*}
    \begin{aligned}\displaystyle
    \limsup_{n \rightarrow +\infty} \gamma \| u_n\|^{2}  < (2+\alpha)\pi.
    \end{aligned}
\end{eqnarray*}
Then there exist $t>1$ close to $1$ and $C>0$ such that
\begin{eqnarray*}
    \begin{aligned}\displaystyle
    \int_{\mathbb{R}^{2}}\Big(e^{\gamma |u_{n}|^{2}}-1\Big)^{t}dx<C, \ \mbox{for all} \ n\in \mathbb{N}.
    \end{aligned}
\end{eqnarray*}
\end{lemma}

By \eqref{snona}, fixed $\tau>3 \ \mbox{and} \ q >2$, for any $\varepsilon>0$ and $\gamma>0$, there exists a constant $\kappa_{\varepsilon}>0$, which depends on $q$, $\gamma$, $\varepsilon$, such that
\begin{equation}\label{snonlina}
|F(u)|\leq\varepsilon|u|^{\tau+1}+\kappa_{\varepsilon}|u|^{q}(e^{\gamma |u|^{2}}-1) \ \mbox{for all} \ u\in H^{1}(\mathbb{R}^{2}),
\end{equation}
and
\begin{equation}\label{snonlinb}
|F(u_{n})|\leq\varepsilon|u_{n}|^{\tau+1}+\kappa_{\varepsilon}|u_{n}|^{q}(e^{\gamma |u_{n}|^{2}}-1), \ \mbox{for all} \ n\in \mathbb{N}.
\end{equation}

By \eqref{nona}, fixed $\tau>3 \ \mbox{and} \ q >2$, for any $\varepsilon>0$ and $\gamma>\gamma_{0}$, there exists a constant $\kappa_{\varepsilon}>0$, which depends on $q$, $\gamma$, $\varepsilon$, such that
\begin{equation}\label{nonlina}
|F(u)|\leq\varepsilon|u|^{\tau+1}+\kappa_{\varepsilon}|u|^{q}(e^{\gamma |u|^{2}}-1) \ \mbox{for all} \ u\in H^{1}(\mathbb{R}^{2}),
\end{equation}
and
\begin{equation}\label{nonlinb}
|F(u_{n})|\leq\varepsilon|u_{n}|^{\tau+1}+\kappa_{\varepsilon}|u_{n}|^{q}(e^{\gamma |u_{n}|^{2}}-1), \ \mbox{for all} \ n\in \mathbb{N}.
\end{equation}

By \eqref{snonlina}, \eqref{snonlinb}, \eqref{nonlina} and \eqref{nonlinb}, we obtain the convergence properties of the nonlinearities is related to the exponential function. Thus, the next lemma is important in our argument.

\begin{lemma}\label{imp} Assume that $(f_{1})$-$(f_{3})$ hold. Let $\{u_{n}\}\subset \mathcal{T}_{r}(a)$ with
\begin{equation}\label{ene}
    \begin{aligned}\displaystyle
	\limsup_{n \rightarrow +\infty} \gamma \| u_n\|^{2}  < (2+\alpha)\pi.
    \end{aligned}
\end{equation}
If $u_{n} \rightharpoonup u$ in $H^{1}_{rad}(\mathbb{R}^{2})$ and $u_n(x) \rightarrow u(x)$ a.e. in $\mathbb{R}$, then
\begin{eqnarray*}
    \begin{aligned}\displaystyle
	|u_n|^{q}(e^{\gamma  |u_n(x)|^{2}}-1) \rightarrow |u|^{q}(e^{\gamma |u(x)|^{2}}-1) \,\, \mbox{in} \,\, L^{1}(\mathbb{R}^{2}).
    \end{aligned}
\end{eqnarray*}
\end{lemma}

\begin{proof}
Setting
\begin{eqnarray*}
    \begin{aligned}\displaystyle
	h_n(x)=e^{\gamma |u_n|^{2}}-1.
    \end{aligned}
\end{eqnarray*}
By \eqref{ene}, there exists $t>1$ close to $1$ such that
\begin{eqnarray*}
    \begin{aligned}\displaystyle
    t\gamma \|u_{n}\|^{2}\leq(2+\alpha)\pi.
    \end{aligned}
\end{eqnarray*}
Thus, by Lemma \ref{TMineq}, we know that
\begin{eqnarray*}
    \begin{aligned}\displaystyle
    \int_{\mathbb{R}^{2}}\Big(e^{\gamma |u_{n}|^{2}}-1\Big)^{t}dx&
    \leq \int_{\mathbb{R}^{2}}\Big(e^{t\gamma |u_{n}|^{2}}-1\Big)dx<C,
    \end{aligned}
\end{eqnarray*}
where $C=C(t,a,\gamma)>0$. Then,
\begin{eqnarray*}
    \begin{aligned}\displaystyle
   \int_{\mathbb{R}^{2}} \Big(h_n(x)\Big)^{t}dx=\int_{\mathbb{R}^{2}}\Big(e^{\gamma |u_{n}|^{2}}-1\Big)^{t}dx<C,
    \end{aligned}
\end{eqnarray*}
which implies that
\begin{eqnarray*}
    \begin{aligned}\displaystyle
	h_{n} \in L^{t}(\mathbb{R}^{2}) \quad \mbox{and} \quad \sup_{n \in \mathbb{N}}|h_n|_{t}<+\infty.
    \end{aligned}
\end{eqnarray*}
Therefore, $\{h_{n}\}$ is a bounded sequence in $L^{t}(\mathbb{R}^{2})$. Thus, for some subequence of $\{u_n\}$, still denoted by itself, we obtain that
	\begin{equation} \label{sNewlimit}
		h_n \rightharpoonup h=e^{\gamma |u|^{2}}-1, \,\, \mbox{in} \,\, L^{t}(\mathbb{R}^{2}).
	\end{equation}
Now, we show that
\begin{equation}\label{ssubimp}
    \begin{aligned}\displaystyle
	|u_{n}|^{q}\rightarrow |u|^{q} \ \ \mbox{in} \ \ L^{t^{\prime}}(\mathbb{R}^{2}),
    \end{aligned}
\end{equation}
where $t^{\prime}=\frac{t}{t-1}$. Then, by the embedding $H_{rad}^{1}(\mathbb{R}^{2}) \hookrightarrow L^{qt'}(\mathbb{R}^{2})$ is compact, we have
\begin{eqnarray*}
    \begin{aligned}\displaystyle
	u_{n} \rightarrow u \quad  L^{qt^{\prime}}(\mathbb{R}^{2}).
    \end{aligned}
\end{eqnarray*}
Hence, we get \eqref{ssubimp}. Together \eqref{sNewlimit} with \eqref{ssubimp}, we know
\begin{eqnarray*}
    \begin{aligned}\displaystyle
    |u_n|^{q}(e^{\gamma  |u_n(x)|^{2}}-1) \rightarrow |u|^{q}(e^{\gamma |u(x)|^{2}}-1) \,\, \mbox{in} \,\, L^{1}(\mathbb{R}^{2}).
    \end{aligned}
\end{eqnarray*}
Then, the proof is complete.
\end{proof}
By Lemma \ref{imp}, we have the following two important Corollaries.
\begin{corollary}\label{weaklim}
Assume that $(f_{1})$-$(f_{3})$ hold. Let $\{u_{n}\}\subset \mathcal{T}_{r}(a)$ with
\begin{equation}
    \begin{aligned}\displaystyle
	\limsup_{n \rightarrow +\infty} \gamma \| u_n\|^{2}  < (2+\alpha)\pi.
    \end{aligned}
\end{equation}
	If $u_{n} \rightharpoonup u$ in $H^{1}_{rad}(\mathbb{R}^{2})$ and $u_{n}(x) \rightarrow u(x)$ a.e in $\mathbb{R}$, then
\begin{eqnarray*}
    \begin{aligned}\displaystyle
	\int_{\mathbb{R}^{2}}(I_{\alpha}\ast F(u_{n}))f(u_{n})\phi dx \rightarrow \int_{\mathbb{R}^{2}}(I_{\alpha}\ast F(u))f(u)\phi dx, \ \mbox{as} \ n\rightarrow\infty,
    \end{aligned}
\end{eqnarray*}
for any $\phi\in C_{0}^{\infty}(\mathbb{R}^{2})$.
\end{corollary}

\begin{proof}
By \cite[Lemma 4.1]{ACTT}, we know
\begin{equation}\label{choq}
    \begin{aligned}\displaystyle
    |I_{\alpha}\ast F(u_{n})|_{\infty}\leq C.
    \end{aligned}
\end{equation}
Hence, for any $\phi\in C_{0}^{\infty}(\mathbb{R}^{2})$, we have
\begin{eqnarray*}
    \begin{aligned}\displaystyle
    |(I_{\alpha}\ast F(u_{n}))f(u_{n})\phi|\leq C|f(u_{n})||\phi|\leq \varepsilon|u_{n}|^{\tau}|\phi| +C|u_{n}|^{q-1}|\phi|(e^{\gamma |u_{n}|^{2}}-1).
    \end{aligned}
\end{eqnarray*}
Let $U=supp\phi$. Then, we obtain
\begin{eqnarray*}
    \begin{aligned}\displaystyle
    \int_{U}|u_{n}|^{\tau}|\phi|dx\rightarrow\int_{U}|u|^{\tau}|\phi|dx, \ \mbox{as} \ n\rightarrow\infty,
    \end{aligned}
\end{eqnarray*}
and
\begin{eqnarray*}
    \begin{aligned}\displaystyle
    \int_{U}|u_{n}|^{q-1}|\phi|(e^{\gamma |u_{n}|^{2}}-1)dx\rightarrow\int_{U}|u|^{q-1}|\phi|(e^{\gamma |u|^{2}}-1)dx, \ \mbox{as} \ n\rightarrow\infty.
    \end{aligned}
\end{eqnarray*}
Now, applying a variant of the Lebesgue Dominated Convergence Theorem, we can deduce that
\begin{eqnarray*}
    \begin{aligned}\displaystyle
	\int_{\mathbb{R}^{2}}(I_{\alpha}\ast F(u_{n}))f(u_{n})\phi dx \rightarrow \int_{\mathbb{R}^{2}}(I_{\alpha}\ast F(u))f(u)\phi dx, \ \mbox{as} \ n\rightarrow\infty,
    \end{aligned}
\end{eqnarray*}
which completes the proof.
\end{proof}

\begin{corollary} \label{Conver} Assume that $(f_{1})$-$(f_{3})$ hold. Let $\{u_{n}\}\subset \mathcal{T}_{r}(a)$ with
\begin{equation}\label{eneb}
    \begin{aligned}\displaystyle
	\limsup_{n \rightarrow +\infty} \gamma \| u_n\|^{2}  < (2+\alpha)\pi.
    \end{aligned}
\end{equation}
	If $u_{n} \rightharpoonup u$ in $H^{1}_{rad}(\mathbb{R}^{2})$ and $u_{n}(x) \rightarrow u(x)$ a.e in $\mathbb{R}$, then
\begin{eqnarray*}
    \begin{aligned}\displaystyle
	\int_{\mathbb{R}^{2}}(I_{\alpha}\ast F(u_{n}))F(u_{n})dx \rightarrow \int_{\mathbb{R}^{2}}(I_{\alpha}\ast F(u))F(u)dx,
    \end{aligned}
\end{eqnarray*}
and
\begin{eqnarray*}
    \begin{aligned}\displaystyle
 \int_{\mathbb{R}^{2}}(I_{\alpha}\ast F(u_{n}))f(u_{n})u_{n}\,dx \rightarrow \int_{\mathbb{R}^{2}}(I_{\alpha}\ast F(u))f(u)u\,dx.
 \end{aligned}
\end{eqnarray*}	
\end{corollary}

\begin{proof} From \eqref{choq}, we know
\begin{eqnarray*}
    \begin{aligned}\displaystyle
    |I_{\alpha}\ast F(u_{n})|_{\infty}\leq C.
    \end{aligned}
\end{eqnarray*}
By \eqref{snona} and \eqref{nona}, we have
\begin{eqnarray*}
    \begin{aligned}\displaystyle
	|F(u_{n})|\leq\varepsilon|u_{n}|^{\tau+1}+C|u_{n}|^{q}(e^{\gamma |u_{n}|^{2}}-1)\,\, \text{ for all }\, u_{n} \in H^{1}(\mathbb{R}^{2}),
    \end{aligned}
\end{eqnarray*}
where $\gamma>\gamma_{0}$, $\tau>3$ and $q>2$. Hence, we have
\begin{eqnarray*}
    \begin{aligned}\displaystyle
    |(I_{\alpha}\ast F(u_{n}))F(u_{n})|\leq C|F(u_{n})|\leq\varepsilon C|u_{n}|^{\tau+1}+C|u_{n}|^{q}(e^{\gamma |u_{n}|^{2}}-1).
    \end{aligned}
\end{eqnarray*}
By 	Lemma \ref{imp}, we know
\begin{eqnarray*}
    \begin{aligned}\displaystyle
	\int_{\mathbb{R}^{2}}|u_n|^{q}(e^{\gamma  |u_n(x)|^{2}}-1)dx \rightarrow \int_{\mathbb{R}^{2}} |u|^{q}(e^{\gamma |u(x)|^{2}}-1)dx \,\, \mbox{as} \,\, n\rightarrow\infty.
    \end{aligned}
\end{eqnarray*}
By the compact embedding $H_{rad}^{1}(\mathbb{R}^{2}) \hookrightarrow L^{p}(\mathbb{R}^{2})$, for $p>2$, we have
\begin{eqnarray*}
    \begin{aligned}\displaystyle
	u_{n}\rightarrow u \quad \mbox{in} \ \ L^{p}(\mathbb{R}^2).
    \end{aligned}
\end{eqnarray*}
Now, applying a variant of the Lebesgue Dominated Convergence Theorem, we can deduce that
\begin{eqnarray*}
    \begin{aligned}\displaystyle
	\int_{\mathbb{R}^{2}}(I_{\alpha}\ast F(u_{n}))F(u_{n})dx \rightarrow \int_{\mathbb{R}^{2}}(I_{\alpha}\ast F(u))F(u)dx  \,\, \mbox{as} \,\, n\rightarrow\infty.
    \end{aligned}
\end{eqnarray*}
	A similar argument works to show that
\begin{eqnarray*}
    \begin{aligned}\displaystyle
	\int_{\mathbb{R}^{2}}(I_{\alpha}\ast F(u_{n}))f(u_{n})u_{n}\,dx \rightarrow \int_{\mathbb{R}^{2}}(I_{\alpha}\ast F(u))f(u)u\,dx \,\, \mbox{as} \,\, n\rightarrow\infty,
    \end{aligned}
\end{eqnarray*}	
which completes the proof.
\end{proof}

\begin{lemma}\label{lamd}
Assume that ($f_{1}$)-($f_{3}$) hold. Let $\{u_{n}\}\subset \mathcal{T}_{r}(a)$ and the $(PS)$ sequence $\{u_{n}\}$ of $J$ is bounded in $H^{1}(\mathbb{R}^{2})$. Then, $\{\lambda_{n}\}$ is a bounded sequence with
\begin{eqnarray*}
    \begin{aligned}\displaystyle
    \limsup\limits_{n\rightarrow\infty}\lambda_{n}=\frac{2+\alpha}{2a^{2}}\int_{\mathbb{R}^{2}}(I_{\alpha}\ast F(u_{n}))F(u_{n})dx.
    \end{aligned}
\end{eqnarray*}	
\end{lemma}

\begin{proof}
By \eqref{EQ2}, we know that
\begin{eqnarray*}
    \begin{aligned}\displaystyle
	-\Delta u_{n}+\lambda_{n}u_{n}=(I_{\alpha}\ast F(u_{n}))f(u_{n})\ + o_n(1),
    \end{aligned}
\end{eqnarray*}
which together with the sequence $\{u_{n}\}$ is bounded in $H^{1}(\mathbb{R}^{2})$ imply that $\{\lambda_{n}\}$ is a bounded sequence.
Thus, by $|u_{n}|_{2}^{2}=a^{2}$, we have
\begin{equation}\label{aineq}
    \begin{aligned}\displaystyle
    \lambda_{n}a^{2}=- |\nabla u_{n}|_{2}^{2}+\int_{\mathbb{R}^{2}}(I_{\alpha}\ast F(u_{n}))f(u_{n})u_{n}dx\ + o_n(1).
    \end{aligned}
\end{equation}
The equality \eqref{aineq} together with the limit \eqref{PEQ2} lead to
\begin{eqnarray*}
    \begin{aligned}\displaystyle
        \limsup\limits_{n\rightarrow\infty}\lambda_{n}=\frac{2+\alpha}{2a^{2}}\int_{\mathbb{R}^{2}}(I_{\alpha}\ast F(u_{n}))F(u_{n})dx,
    \end{aligned}
\end{eqnarray*}
which completes the proof.
\end{proof}

\begin{lemma}\label{uni}
Assume that $(f_{1})$-$(f_{3})$ and $(f_{5})$ hold with $u\in \mathcal{T}_{r}(a)$. Then the function $\tilde{J}_{u}(s)=J(\mathcal{H}(u,s))$ reaches its unique maximum at a point $s(u)\in \mathbb{R}$ such that $\mathcal{H}(u,s(u))\in \mathcal{P}(a)$.
\end{lemma}
\begin{proof}
For $u\in \mathcal{T}_{r}(a)$ and $s\in \mathbb{R}$, we know
\begin{eqnarray*}
    \begin{aligned}\displaystyle
    \tilde{J}_{u}(s)&=\frac{1}{2}\int_{\mathbb{R}^{2}}|\nabla \mathcal{H}(u,s)|^{2}dx-\frac{1}{2}\int_{\mathbb{R}^{2}} (I_{\alpha}\ast F(\mathcal{H}(u,s)))F(\mathcal{H}(u,s)) dx\\
    &=\frac{e^{2s}}{2}\int_{\mathbb{R}^{2}}|\nabla u|^{2}dx-\frac{1}{2e^{(2+\alpha)s}}\int_{\mathbb{R}^{2}} (I_{\alpha}\ast F(e^{s}u))F(e^{s}u) dx.
    \end{aligned}
\end{eqnarray*}
Then, we have
\begin{equation}\label{DS1}
    \begin{aligned}\displaystyle
\tilde{J}^{\prime}_{u}(s)=&e^{2s}\int_{\mathbb{R}^{2}}|\nabla u|^{2}dx-\frac{2+\alpha}{2e^{(2+\alpha)s}}\int_{\mathbb{R}^{2}} (I_{\alpha}\ast F(e^{s}u))F(e^{s}u) dx\\
&-\frac{2}{2e^{(2+\alpha)s}}\int_{\mathbb{R}^{2}} I_{\alpha}\ast F(e^{s}u))f(e^{s}u)e^{s}u dx\\
=& e^{2s}(\int_{\mathbb{R}^{2}}|\nabla u|^{2}dx-\psi(s)),
    \end{aligned}
\end{equation}
where
\begin{eqnarray*}
    \begin{aligned}\displaystyle
    \psi(s)=\int_{\mathbb{R}^{2}} (_{\alpha}\ast\frac{F(e^{s}u)}{(e^{s})^{2+\frac{\alpha}{2}}})\frac{\tilde{F}(e^{s}u)}{(e^{s})^{2+\frac{\alpha}{2}}}.
    \end{aligned}
\end{eqnarray*}
By Lemmas \ref{geo} and \ref{P1}, we know that there exists at least a $s_{0}\in \mathbb{R}$ such that $f^{\prime}_{u}(s)_{|_{s=s_{0}}}=0$. For any $t\in \mathbb{R}$, $t\neq0$, from $(f_{3})$ and $(f_{5})$, we see that $\frac{F(st)}{s^{2+\frac{\alpha}{2}}}$ is strictly increasing in $s\in(0,\infty)$ and $\frac{\tilde{F}(st)}{s^{2+\frac{\alpha}{2}}}$ is nondecreasing in $s\in(0,\infty)$. This implies $ \psi(s)$ is strictly increasing in $s\in(0,\infty)$ and there is at most one $s(u)\in \mathbb{R}$ such that $\mathcal{H}(u,s(u))\in \mathcal{P}(a)$. Thus, for $u\in \mathcal{T}_{r}(a)$, there exists a unique maximum at a point $s(u)\in \mathbb{R}$ such that $\mathcal{H}(u,s(u))\in \mathcal{P}(a)$.
\end{proof}

\begin{lemma}\label{fin}
Assume that $(f_{1})$-$(f_{3})$ and $(f_{5})$ hold with $u\in \mathcal{T}_{r}(a)$. Then
\begin{eqnarray*}
    \begin{aligned}\displaystyle
    m(a)= c(a):=\inf\limits_{u\in \mathcal{P}(a)}\max\limits_{s\in \mathbb{R}}J(\mathcal{H}(u,s)).
    \end{aligned}
\end{eqnarray*}
\end{lemma}
\begin{proof}
For $u\in \mathcal{T}_{r}(a)$, by Lemma \ref{uni}, we know that
\begin{eqnarray*}
    \begin{aligned}\displaystyle
    \max\limits_{s\in \mathbb{R}}J(\mathcal{H}(u,s))=J(\mathcal{H}(u,s(u)))\geq m(a).
    \end{aligned}
\end{eqnarray*}
Thus, $m(a)\leq c(a)$. On the other hand, for $u\in \mathcal{P}(a)$, Lemma \ref{uni} also implies that
\begin{eqnarray*}
    \begin{aligned}\displaystyle
J(u)=\max\limits_{s\in \mathbb{R}}J(\mathcal{H}(u,s))\geq \inf\limits_{u\in \mathcal{P}(a)}\max\limits_{s\in \mathbb{R}}J(\mathcal{H}(u,s)).
    \end{aligned}
\end{eqnarray*}
Hence, $m(a)= c(a)$.
\end{proof}

\section{{\bfseries Proof of Theorem \ref{T1}}}\label{PFTa}
In this section, we assume that $f$ has subcritical growth and restrict our study in $H_{rad}^{1}(\mathbb{R}^{2})$.\\
{\bf Proof of Theorem} $\ref{T1}$. First we show $u_{n}\rightharpoonup u$ in $H_{rad}^{1}(\mathbb{R}^{2})$, where $u\neq 0$. By Lemma \ref{bounded}, we have $u_{n}\rightharpoonup u$ in $H_{rad}^{1}(\mathbb{R}^{2})$. For $\gamma > 0$ close to $0$ and the sequence $\{u_{n}\}$ is bounded, we obtain
\begin{eqnarray*}
    \begin{aligned}\displaystyle
	\limsup_{n \rightarrow +\infty} \gamma \| u_n\|^{2}  < (2+\alpha)\pi.
    \end{aligned}
\end{eqnarray*}
Then,  using Corollary \ref{Conver}, it follows that
\begin{equation}\label{limitaa}
    \begin{aligned}\displaystyle
\lim_{n \rightarrow +\infty}\int_{\mathbb{R}^2}(I_{\alpha}\ast F(u_{n}))f(u_{n})u_{n}\,dx=\int_{\mathbb{R}^2}(I_{\alpha}\ast F(u))f(u)u\,dx,
    \end{aligned}
\end{equation}
and
\begin{equation}\label{limitab}
    \begin{aligned}\displaystyle
\lim_{n \rightarrow +\infty}\int_{\mathbb{R}^2}(I_{\alpha}\ast F(u_{n}))F(u_{n})\,dx=\int_{\mathbb{R}^2}(I_{\alpha}\ast F(u))F(u)\,dx,
    \end{aligned}
\end{equation}
where $u_{n} \rightharpoonup u$ in $H_{rad}^{1}(\mathbb{R}^{2})$. The last limit implies that $u \not=0$, because otherwise, Corollary \ref{Conver} gives
\begin{eqnarray*}
    \begin{aligned}\displaystyle
\lim_{n \rightarrow +\infty}\int_{\mathbb{R}^2}(I_{\alpha}\ast F(u_{n}))F(u_{n})\,dx=\lim_{n \rightarrow +\infty}\int_{\mathbb{R}^2}(I_{\alpha}\ast F(u_{n}))f(u_{n})u_{n}\,dx=0,
    \end{aligned}
\end{eqnarray*}
and by Lemma \ref{lamd},
\begin{eqnarray*}
    \begin{aligned}\displaystyle
\limsup\limits_{n\rightarrow\infty}\lambda_{n}=\frac{2+\alpha}{2a^{2}}\int_{\mathbb{R}^{2}}(I_{\alpha}\ast F(u_{n}))F(u_{n})dx>0.
    \end{aligned}
\end{eqnarray*}
Since $\{u_{n}\}$ is bounded in $H_{rad}^{1}(\mathbb{R}^{2})$, $\lambda_{n}>0$, Corollary \ref{Conver} and the equality below
\begin{equation*}
	|\nabla u_n|^{2}_{2}+\lambda_{n}| u_{n}|_{2}^{2}=\int_{\mathbb{R}^2}(I_{\alpha}\ast F(u_{n}))f(u_{n})u_{n}\,dx+o_n(1),
\end{equation*}
lead to
\begin{equation}  \label{m2}
	-\lambda_{n}a^{2}=|\nabla u_{n}|_{2}^{2}+o_n(1).
\end{equation}
From this,
\begin{eqnarray*}
    \begin{aligned}\displaystyle
0 \geq \limsup_{n \rightarrow +\infty} (-\lambda_{n})a^{2}= \limsup_{n \rightarrow +\infty} |\nabla u_{n}|^{2}_{2} \geq \liminf_{n \rightarrow +\infty} |\nabla u_{n}|^{2}_{2}\geq 0,
    \end{aligned}
\end{eqnarray*}
thus $|\nabla u_n|^{2}_{2} \to 0$, which is absurd, because $m(a)>0$.

Then, we show that $\lambda>0$. By Lemma \ref{lamd} and $(f_{3})$, there exists a bounded sequence $\{\lambda_{n}\}$ such that
\begin{eqnarray*}
    \begin{aligned}\displaystyle
\limsup\limits_{n\rightarrow\infty}\lambda_{n}=\frac{2+\alpha}{2a^{2}}\int_{\mathbb{R}^{2}}(I_{\alpha}\ast F(u_{n}))F(u_{n})dx>0.
    \end{aligned}
\end{eqnarray*}
From this for some subsequence, still denoted by $\{\lambda_{n}\}$, we can assume that
\begin{eqnarray*}
    \begin{aligned}\displaystyle
    \lambda_{n}\rightarrow\lambda>0 \ \mbox{as} \ n\rightarrow\infty.
    \end{aligned}
\end{eqnarray*}
Therefore $\lambda>0$. Then, by Corollary \ref{weaklim}, the equality \eqref{EQ2} implies that
\begin{eqnarray*}
    \begin{aligned}\displaystyle
-\Delta u+\lambda u=(I_{\alpha}\ast F(u))f(u) \ \ \mbox{in} \ \mathbb{R}^{2}.
    \end{aligned}
\end{eqnarray*}
Thus, we deduce that $P(u)=0$.
Now, we obtain that $u_{n}\rightharpoonup u\neq0$. Then we show the strong convergence that $u_{n}\rightarrow u$ in $H_{rad}^{1}(\mathbb{R}^{2})$. The proof is divided into two steps.

{\bf Step 1.} We show that $\lim\limits_{n\rightarrow\infty} |\nabla u_{n}|_{2}^{2}=|\nabla u|_{2}^{2}$.

By \eqref{Pohozaev}, \eqref{limitaa} and \eqref{limitab}, together with the weak convergence
\begin{eqnarray*}
    \begin{aligned}\displaystyle
 P(u_{n})-P(u)=o_{n}(1),
    \end{aligned}
\end{eqnarray*}
we deduce
\begin{eqnarray*}
    \begin{aligned}\displaystyle
\lim\limits_{n\rightarrow\infty}|\nabla{u}_{n}|_{2}^{2}-|\nabla{u}|_{2}^{2}=0,
    \end{aligned}
\end{eqnarray*}
that is
\begin{eqnarray*}
    \begin{aligned}\displaystyle
\lim\limits_{n\rightarrow\infty}|\nabla u_{n}|_{2}^{2}=|\nabla u|_{2}^{2} .
    \end{aligned}
\end{eqnarray*}

{\bf Step 2.} We show that $|u|_{2}=a$.

Combining \eqref{mu2}, \eqref{EQ2}, \eqref{PEQ2} and $P(u)=0$, we obtain
\begin{eqnarray*}
    \begin{aligned}\displaystyle
    \lambda a^{2}=\lambda|u|_{2}^{2}.
    \end{aligned}
\end{eqnarray*}
Then, $0<|u|_{2}\leq a$ implying that $|u|_{2}=a$. Thus, $u_{n}\rightarrow u$ in $H_{rad}^{1}(\mathbb{R}^{2})$. Finally, by Lemma \ref{fin}, we obtain that $u$ is a normalized ground state solution of \eqref{aa}.

\section{{\bfseries On the mini-max level}}\label{minimax}
In this section, we obtain an upper bound for the minimax level. Thus, we obtain an upper bound for $\|u\|^{2}$, which is important for exponential critical problem.

\begin{lemma} \label{ESTMOUNTPASS} Assume that ($f_{4}$) holds, there holds $\displaystyle \lim_{\mu \rightarrow +\infty}m(a)=0$.	
\end{lemma}

\begin{proof} In what follow we set the path $h_{0}(t)=\mathcal{H}\big(u_{0}, (1-t)s_1+ts_2\big) \in \Gamma$. Then, by $(f_{4})$,
\begin{eqnarray*}
    \begin{aligned}\displaystyle
	m(a) \leq \max_{t \in [0,1]}J(h_0(t)) &\leq \max_{s \in \mathbb{R}}\left\{\frac{e^{2s}}{2}|\nabla u_{0}|_{2}^{2}-\frac{C\mu}{\sigma}e^{(2\sigma-2-\alpha)s}|u_{0}|_{\frac{8\sigma}{2+\alpha}}^{4\sigma}\right\}\\
&\leq \max_{r\geq0}\left\{\frac{r^{2}}{2}|\nabla u_{0}|_{2}^{2}-\frac{C\mu}{\sigma}r^{2\sigma-2-\alpha}|u_{0}|_{\frac{8\sigma}{2+\alpha}}^{4\sigma}\right\}.
    \end{aligned}
\end{eqnarray*}
Thus,
\begin{eqnarray*}
    \begin{aligned}\displaystyle
	m(a) \leq C_{1}\left(\frac{1}{\mu}\right)^{\frac{2}{2\sigma-4-\alpha}} \to 0 \quad \mbox{as} \quad \mu \to +\infty,
    \end{aligned}
\end{eqnarray*}
for some $C_1>0$. Here, we have used the fact that $\sigma>2+\frac{\alpha}{2}$.
\end{proof}

\begin{lemma}\label{newlem1} Assume that ($f_{4}$) holds. Let $\{u_{n}\}\subset \mathcal{T}_{r}(a)$ is the $(PS)$ sequence of $J$. Then there holds
\begin{eqnarray*}
    \begin{aligned}\displaystyle
	\limsup_{n \rightarrow +\infty}(I_{\alpha}\ast F(u_{n}))F(u_{n})\,dx \leq \frac{4-2\alpha}{2\theta-4-\alpha}m(a).
    \end{aligned}
\end{eqnarray*}
\end{lemma}
\begin{proof} Using the fact that $J(u_{n})=m(a)+o_n(1)$ and $P(u_{n})=o_n(1)$, it follows that
\begin{eqnarray*}
    \begin{aligned}\displaystyle
	(2+\alpha){J}(u_{n})+P(u_{n})=(2+\alpha)m(a)+o_n(1),
    \end{aligned}
\end{eqnarray*}
and so,
\begin{equation}\label{equa1}
    \begin{aligned}\displaystyle
	\frac{4+\alpha}{2}|\nabla u_{n}|_{2}^{2}-\int_{\mathbb{R}^2} (I_{\alpha}\ast F(u_{n}))f(u_{n})u_{n}\,dx=(2+\alpha)m(a)+o_n(1).
    \end{aligned}
\end{equation}
Using that $J(u_{n})=m(a)+o_n(1)$, we obtain
\begin{equation}\label{equa2}
    \begin{aligned}\displaystyle
|\nabla u_{n}|_{2}^{2}=\int_{\mathbb{R}^2} (I_{\alpha}\ast F(u_{n}))F(u_{n}) dx+2m(a)+o_n(1).
    \end{aligned}
\end{equation}
Together \eqref{equa1} and \eqref{equa2}, we get
\begin{eqnarray*}
    \begin{aligned}\displaystyle
(4+\alpha)\int_{\mathbb{R}^2} (I_{\alpha}\ast F(u_{n}))F(u_{n}) dx+4m(a)+o_n(1)-2\int_{\mathbb{R}^2}(I_{\alpha}\ast F(u_{n}))f(u_{n})u_{n}\,dx=o_n(1).
    \end{aligned}
\end{eqnarray*}
Then, by $(f_{3})$,
\begin{eqnarray*}
    \begin{aligned}\displaystyle
4m(a)+o_n(1)&=2\int_{\mathbb{R}^2}(I_{\alpha}\ast F(u_{n}))f(u_{n})u_{n}\,dx-(4+\alpha)\int_{\mathbb{R}^2} (I_{\alpha}\ast F(u_{n}))F(u_{n})  dx\\
&\geq (2\theta-4-\alpha)\int_{\mathbb{R}^2} (I_{\alpha}\ast F(u_{n}))F(u_{n}) dx.
    \end{aligned}
\end{eqnarray*}
Since $\theta>2+\frac{\alpha}{2}$, we have
\begin{eqnarray*}
    \begin{aligned}\displaystyle
\limsup_{n \rightarrow +\infty}\int_{\mathbb{R}^2}(I_{\alpha}\ast F(u_{n}))F(u_{n})\,dx \leq \frac{4}{2\theta-4-\alpha}m(a).
    \end{aligned}
\end{eqnarray*}
Then, the proof is complete.
\end{proof}

\begin{lemma} \label{boundd}Assume that ($f_{4}$) holds. Let $\{u_{n}\}\subset \mathcal{T}_{r}(a)$. Then $\{u_{n}\}$ is the $(PS)$ sequence of $J$ satisfies
\begin{eqnarray*}
    \begin{aligned}\displaystyle
\displaystyle \limsup_{n \rightarrow +\infty}|\nabla u_{n}|_{2}^{2} \leq \frac{4\theta-4-2\alpha}{2\theta-4-\alpha} m(a).
    \end{aligned}
\end{eqnarray*}
Hence, there exists $\mu^*>0$ such that
\begin{equation}\label{level}
    \begin{aligned}\displaystyle
	\limsup_{n \to +\infty}|\nabla u_{n}|_{2}^{2}<\frac{(2+\alpha)\pi}{\gamma_{0}}-a^{2}, \,\,\text{for any} \,\,\,\mu \geq \mu^*.
    \end{aligned}
\end{equation}
\end{lemma}

\begin{proof}  Since $J(u_{n})=m(a)+o_n(1)$, we have
\begin{eqnarray*}
    \begin{aligned}\displaystyle
	\int_{\mathbb{R}^2}|\nabla u_{n}|^2\,dx=2m(a)+\int_{\mathbb{R}^2}(I_{\alpha}\ast F(u_{n}))F(u_{n})\,dx+o_n(1).
    \end{aligned}
\end{eqnarray*}
Thus, by Lemma \ref{newlem1}, we obtain
\begin{eqnarray*}
    \begin{aligned}\displaystyle
	\limsup_{n \to +\infty}|\nabla u_n|_{2}^{2}  \leq \frac{4\theta-4-2\alpha}{2\theta-4-\alpha} m(a),
    \end{aligned}
\end{eqnarray*}
which implies \eqref{level} holds.
\end{proof}

\section{{\bfseries Proof of Theorem \ref{T2}}}\label{PFTb}
In this section, we assume that $f$ has critical growth and restrict our study in $H_{rad}^{1}(\mathbb{R}^{2})$.\\
{\bf Proof of Theorem} $\ref{T2}$.
By Lemma \ref{boundd}, we have
\begin{eqnarray*}
    \begin{aligned}\displaystyle
 \limsup_{n \to +\infty}|\nabla u_{n}|_{2}^{2}<\frac{(2+\alpha)\pi}{\gamma_{0}}-a^{2}, \,\,\text{for any} \,\,\,\mu \geq \mu^*,
    \end{aligned}
\end{eqnarray*}
which implies that
\begin{eqnarray*}
    \begin{aligned}\displaystyle
    \limsup_{n \to +\infty}\| u_{n}\|^{2}<\frac{(2+\alpha)\pi}{\gamma_{0}}, \,\,\text{for any} \,\,\,\mu \geq \mu^*.
    \end{aligned}
\end{eqnarray*}
Hence, for $\gamma>\gamma_{0}$ close to $\gamma_{0}$, we obtain
\begin{eqnarray*}
    \begin{aligned}\displaystyle
	\limsup_{n \rightarrow +\infty} \gamma \| u_n\|^{2}  < (2+\alpha)\pi.
    \end{aligned}
\end{eqnarray*}
Following a similar argument as Section \ref{PFTa}, we complete the proof of Theorem \ref{T2}.
\qed

\smallskip

\smallskip
\smallskip

\noindent{\bfseries Acknowledgements:}
The research has been supported by National Natural Science Foundation of China 11971392, Natural Science
Foundation of Chongqing, China cstc2019jcyjjqX0022.

\end{document}